\documentclass[a4paper,11pt]{amsart}
\usepackage{graphicx}
\usepackage{mathptmx}
\usepackage{amsmath}
\usepackage{amssymb}
\usepackage{enumitem}
\usepackage{xcolor}
\usepackage{pgfplots}

\newmuskip\pFqmuskip

\newcommand*\pFq[6][8]{%
  \begingroup % only local assignments
  \pFqmuskip=#1mu\relax
  \mathcode`=\string"8000
  \begingroup\lccode`\~=`\,
  \lowercase{\endgroup\let~}\pFqcomma
  F^{#2}_{#3}{\left(\genfrac..{0pt}{}{#4}{#5}\bigg|#6\right)}%
  \endgroup
}
\newcommand{\pFqcomma}{\mskip\pFqmuskip}

\newtheorem{theorem}{Theorem}[section]

\newtheorem{corollary}[theorem]{Corollary}

\begin{document}

\title[On Generalized Chebyshev polynomials]{On Generalized Chebyshev polynomials}

\author{Taekyun  Kim}
\address{Department of Mathematics, Kwangwoon University, Seoul 139-701, Republic of Korea}
\email{tkkim@kw.ac.kr}
\author{Dae San  Kim}
\address{Department of Mathematics, Sogang University, Seoul 121-742, Republic of Korea}
\email{dskim@sogang.ac.kr}


\maketitle

\begin{abstract}
In this paper, we introduce and study a novel family of generalized Chebyshev polynomials $K_{n}(x|a,b,c)$ defined via a rational generating function. We demonstrate that the classical Chebyshev polynomials of the first, second, third, and fourth kinds naturally emerge as special cases of this framework. Furthermore, we derive comprehensive recurrence relations, tridiagonal determinant representations, and explicit closed-form expressions for these polynomials. We also establish some connections between the generalized Chebyshev polynomials and other classical sequences, including Morgan-Voyce polynomials, Fibonacci polynomials, and Fubini polynomials via Stirling numbers of the second kind. Finally, we examine the associated Euler-Seidel matrix and provide a continued fraction expansion for the quotient of consecutive polynomial terms.
\end{abstract}

\section{Introduction}
Polynomial sequences and their generalizations play a fundamental role in various branches of mathematics, including combinatorics, number theory, and mathematical analysis \cite{3}. Among these, classic objects of study with rich algebraic and combinatorial properties include the Chebyshev polynomials of the four kinds \cite{1, 15}, Fibonacci numbers and polynomials \cite{2, 12}, Fubini polynomials \cite{16, 19}, and Morgan-Voyce polynomials \cite{13, 17, 18}. Recent advances in combinatorial sequences and special functions have further highlighted the utility of related structures, such as special number systems \cite{4, 6, 7, 8}, interconnections between Chebyshev and Fibonacci families \cite{9, 11}, geometric connections to polynomials \cite{14}, and systematic matrix methods like the Euler-Seidel approach \cite{5, 10}. \par
In this paper, we introduce and investigate a unified framework for these polynomial families by defining the generalized Chebyshev polynomials $K_{n}(x\vert{}a,b,c)$ through the generating function:
\begin{equation*}
\frac{1+cxt}{1+axt+bt^{2}}=\sum_{n=0}^{\infty}K_{n}(x\vert{}a,b,c)t^{n},
\end{equation*}
where $a,b,c\in\mathbb{R}$ (or $\mathbb{C}$). We also examine the specialized subfamily $T_{n}(x\vert{}a,b) = K_{n}(x\vert{}a,b,0)$. A notable feature of this generalization is that classical Chebyshev polynomials emerge as natural special cases under specific parameter choices \cite{1}:
\begin{align*}
&K_{n}\big(x\big\vert{}-2,1,-1\big)=T_{n}(x),\quad K_{n}\big(x\big\vert{}-2,1,0\big)=U_{n}(x), \\
&K_{n}\bigg(x\bigg\vert{}-2,1,-\frac{1}{x}\bigg)=V_{n}(x),\quad K_{n}\bigg(x\vert{}-2,1,\frac{1}{x}\bigg)=W_{n}(x).
\end{align*} \par
The primary contributions of this work encompass structural properties, determinant representations, and connections to other well-known sequences. Specifically, we derive recurrence relations and establish a tridiagonal determinant representation for $K_{n}(x\vert{}a,b,c)$ (Theorem 2.1). Furthermore, we establish closed-form expressions for both $T_{n}(x\vert{}a,b)$ (Theorem 2.3) and the broader family $K_{n}(x\vert{}a,b,c)$ (Theorem 2.4). We also connect $T_{n}(x\vert{}a,b)$ to Morgan-Voyce polynomials \cite{13, 17, 18}, demonstrating that $T_{2n}(x\vert{}1,-1)=b_{n}(x^{2})$ and $T_{2n+2}\big(x\big\vert{}1,-1\big)-T_{2n}\big(x\big\vert{}1,-1\big)=x^{2}B_{n}(x^{2})$ (Theorem 2.5), and link specialized parameters to Fibonacci polynomials $f_{n}(x)$ \cite{9, 11, 12} (Corollary 2.7). Finally, we express Fubini polynomials $F_{n}(x)$ \cite{16, 19} in terms of the generalized Chebyshev polynomials and Stirling numbers of the second kind \cite{3, 7} (Theorem 2.8).This paper is structured as follows:\begin{itemize}\item Section 1 recalls the background material, including the Fibonacci numbers $f_{n}$ and polynomials $f_{n}(x)$ \cite{2, 12}, Stirling numbers of the first and second kinds \cite{3}, Fubini polynomials $F_{n}(x)$ \cite{16, 19}, and the four kinds of Chebyshev polynomials along with their closed-form expressions and tridiagonal determinant representations \cite{1, 15}.\item Section 2 contains the main results, detailing the definition, recurrence relations, determinant representations, and closed-form expressions for the generalized Chebyshev polynomials $K_{n}(x\vert{}a,b,c)$ and $T_{n}(x\vert{}a,b)$, alongside their connections to Morgan-Voyce polynomials \cite{13, 14, 17, 18}, Fibonacci polynomials \cite{9, 11}, and Fubini polynomials \cite{16, 19}.\item Section 3 investigates the Euler-Seidel matrix \cite{5, 10} associated with the initial sequence $\left(K_{n}(x\vert{}a,b,c)\right)_{n\ge 0}$ and presents a continued fraction expression for the quotient ratio $\frac{K_{n}(x\vert{}a,b,c)}{K_{n-1}(x\vert{}a,b,c)}$.\end{itemize} \par
The Fibonacci numbers are defined by \cite{2, 12}
\begin{equation}
f_{0}=0,\ f_{1}=1\quad \mathrm{and}\quad f_{n+2}=f_{n+1}+f_{n},\quad (n\ge 0). \label{1}
\end{equation}
From \eqref{1}, we have
\begin{equation*}
f_{n}=\frac{1}{\beta-\alpha}\big(\beta^{n}-\alpha^{n}\big)=\sum_{k=0}^{[\frac{n-1}{2}]}\binom{n-1-k}{k},\quad (n\ge 1),
\end{equation*}
where
\begin{equation*}
\beta=\frac{1+\sqrt{5}}{2},\quad \alpha=\frac{1-\sqrt{5}}{2}.
\end{equation*}
The Fibonacci polynomials are given by
\begin{equation}
f_{0}(x)=0,\ f_{1}(x)=1,\quad f_{n+2}(x)=xf_{n+1}(x)+f_{n}(x),\ (n\ge 0). \label{2}
\end{equation}
Thus, by \eqref{2}, we get \cite{2, 12}
\begin{equation}
\frac{t}{1-xt-t^{2}}=\sum_{n=0}^{\infty}f_{n}(x)t^{n}.\label{3}
\end{equation}
From \eqref{3}, we note that
\begin{equation}
f_{n}(x)=\sum_{k=0}^{[\frac{n-1}{2}]}\binom{n-k-1}{k}x^{n-2k-1},\ (n\ge 1).\label{4}
\end{equation} \par
The Stirling numbers of the first kind are defined as \cite{3},
\begin{equation}
(x)_{n}=\sum_{k=0}^{n}S_{1}(n,k)x^{k},\,\, (n\ge 0), \quad \frac{1}{k!}\big(\log(1+t)\big)^{k}=\sum_{n=k}^{\infty}S_{1}(n,k)\frac{t^{n}}{n!}, \,\, (k \ge 0), \label{5}
\end{equation}
where
\begin{equation*}
(x)_{0}=1,\quad (x)_{n}=x(x-1)\cdots (x-n+1),\ (n\ge 1).
\end{equation*}
As the inversion formula of \eqref{5}, the Stirling numbers of the second kind are given by \cite{3}
\begin{equation}
x^{n}=\sum_{k=0}^{n}{n \brace k}(x)_{k},\,\, (n\ge 0),\quad \frac{1}{k!}(e^{t}-1)^{k}=\sum_{n=k}^{\infty}S_{2}(n,k)\frac{t^{n}}{n!},\,\,(k \ge 0). \label{6}
\end{equation}
The Fubini polynomials (or geometric polynomials) are defined by \cite{16, 19}
\begin{equation}
F_{n}(x)=\sum_{k=0}^{n}{n \brace k}k!x^{k},\quad (n\ge 0). \label{7}
\end{equation}
When $x=1,\ F_{n}=F_{n}(1)$ are called Fubini numbers. \par
The Chebyshev polynomials of the first kind $T_{n}(x)$ and of second kind $U_{n}(x)$ are respectively  given by \cite{1, 15}
\begin{equation}
T_{0}(x)=1,\ T_{1}(x)=x,\quad T_{n+1}(x)=2xT_{n}(x)-T_{n-1}(x),\ (n\ge 1),\label{8}	
\end{equation}
and
\begin{equation}
U_{0}(x)=1,\ U_{1}(x)=2x,\quad U_{n+1}(x)=2xU_{n}(x)-U_{n-1}(x),\ (n\ge 1).\label{9}
\end{equation}
Thus, by \eqref{9}, we get
\begin{equation}
U_{n}(x)=\frac{1}{2\sqrt{x^{2}-1}}\Big((x+\sqrt{x^{2}-1})^{n+1}-(x-\sqrt{x^{2}-1})^{n+1}\Big),\label{10}
\end{equation}
where $n$ is a nonnegative integer. \\
It is easy to show that
\begin{equation*}
i^{-n}U_{n}\bigg(\frac{i}{2}\bigg)=F_{n+1},\ (n\ge 0),\quad \mathrm{where}\quad i=\sqrt{-1}.
\end{equation*}
Let $V_{n}(x)$ and $W_{n}(x)$ be respectively the Cheyshev polynomials of the third and fourth kinds. Then we have \cite{1, 15}
\begin{equation}
V_{0}(x)=1,\ V_{1}(x)=2x-1,\quad V_{n}(x)=2xV_{n-1}(x)-V_{n-2}(x),\ (n\ge 2),\label{11}
\end{equation}
and
\begin{equation}
W_{0}(x)=1,\ W_{1}(x)=2x+1,
W_{n}(x)=2xW_{n-1}(x)-W_{n-2}(x),\ (n\ge 2). \label{12}
\end{equation}
By \eqref{11} and \eqref{12}, $V_{n}(x)$ and $W_{n}(x)$ can be expressed by the following $n \times n$ tridiagonal determinants:
\begin{equation*}
V_{n}(x)=\left|\begin{matrix}
2x & 1 & 0 & 0 & \dots & 0  \\
1 & 2x & 1 & 0 & \dots & 0  \\
0 & 1 & 2x & 1 & \dots & 0  \\
\vdots & \vdots &  \ddots & \ddots & \ddots & \vdots \\
0 & \dots& 0 & 1 & 2x & 1\\
0 & \dots & 0 & 0 & 1 & 2x-1
\end{matrix}\right|,
\end{equation*}
and
\begin{equation*}
W_{n}(x)=\left|\begin{matrix}
2x & 1 & 0 & 0 & \dots & 0  \\
1 & 2x & 1 & 0 & \dots & 0  \\
0 & 1 & 2x & 1 & \dots & 0  \\
\vdots & \vdots &  \ddots & \ddots & \ddots & \vdots \\
0 & \dots& 0 & 1 & 2x & 1\\
0 & \dots & 0 & 0 & 1 & 2x+1
\end{matrix}\right|.
\end{equation*}

\section{Generalized Chebyshev polynomials}
For $a,b,c\in\mathbb{R}$ (or $\mathbb{C}$), let us consider the generalized Chebyshev polynomials given by
\begin{equation}
\frac{1+cxt}{1+axt+bt^{2}}=\sum_{n=0}^{\infty}K_{n}(x|a,b,c)t^{n}. \label{13}	
\end{equation}
Then, by \eqref{13}, we get
\begin{equation}
\begin{aligned}
1+cxt&=K_{0}\big(x|a,b,c\big)+\Big(K_{1}\big(x|a,b,c\big)+axK_{0}\big(x|a,b,c\big)\Big)t \\
&\quad +\sum_{n=2}^{\infty}\Big(K_{n}\big(x|a,b,c\big)+axK_{n-1}\big(x|a,b,c\big)+bK_{n-2}\big(x|a,b,c\big)\Big)t^{n}.
\end{aligned}\label{14}
\end{equation}
Comparing the coefficients on both sides of \eqref{14}, we have
\begin{equation}
K_{0}\big(x|a,b,c\big)=1,\quad K_{1}\big(x|a,b,c\big)=-ax+cx,	\label{15}
\end{equation}
and
\begin{equation}
K_{n}\big(x|a,b,c\big)=-axK_{n-1}\big(x|a,b,c\big)-bK_{n-2}\big(x|a,b,c\big),\ (n\ge 2).\label{16}
\end{equation}
From \eqref{15} and \eqref{16}, we obtain the next theorem.
\begin{theorem}
For $n\in\mathbb{N}$, $K_{n}\big(x|a,b,c\big)$ can be computed as the determinant of the following $n \times n$ tridiagonal matrix.
\begin{equation*}
K_{n}\big(x|a,b,c\big)=\left|\begin{matrix}
-ax & b & 0 & 0 & \dots & 0  \\
1 & -ax & b & 0 & \dots & 0  \\
0 & 1 & -ax & b & \dots & 0  \\
\vdots & \vdots &  \ddots & \ddots & \ddots & \vdots \\
0 &\dots& 0 & 1 & -ax &  b \\
0 & \dots & 0 & 0 & 1 &  -ax+cx
\end{matrix}\right|.
\end{equation*}
\end{theorem}
\emph{Remark.} Note from \eqref{8}, \eqref{9}, \eqref{11} and \eqref{12} that
\begin{align*}
&K_{n}\big(x\big|-2,1,-1\big)=T_{n}(x),\ K_{n}\big(x\big|-2,1,0\big)=U_{n}(x), \\
&K_{n}\bigg(x\bigg|-2,1,-\frac{1}{x}\bigg)=V_{n}(x),\ K_{n}\bigg(x|-2,1,\frac{1}{x}\bigg)=W_{n}(x).
\end{align*} \par
Now, we define the generalized Chebyshev polynomials of the second kind by
\begin{equation}
\frac{1}{1+axt+bt^{2}}=\sum_{n=0}^{\infty}T_{n}\big(x|a,b\big)t^{n}. \label{17}	
\end{equation}
Thus, by \eqref{17}, we get
\begin{equation}
\begin{aligned}
&T_{0}\big(x|a,b\big)=1,\quad T_{1}\big(x\big|a,b\big)=-ax, \\
&T_{n}\big(x\big|a,b\big)=-axT_{n-1}\big(x\big|a,b\big)-bT_{n-2}\big(x|a,b\big),\ (n\ge 2).
\end{aligned}\label{18}
\end{equation}
From \eqref{18} or Theorem 2.1, we have the following result.
\begin{corollary}
For $n\in\mathbb{N}$, $T_{n}\big(x\big|a,b\big)$ is given by the determinant of the following $n \times n$ tridiagonal matrix.
\begin{equation*}
T_{n}\big(x|a,b\big)=\left|\begin{matrix}
-ax & b & 0 & 0 & \dots & 0  \\
1 & -ax & b & 0 & \dots & 0  \\
0 & 1 & -ax & b & \dots & 0  \\
\vdots & \vdots &  \ddots & \ddots & \ddots & \vdots \\
0 &\dots& 0 & 1 & -ax &  b \\
0 & \dots & 0 & 0 & 1 & -ax
\end{matrix}\right|.
\end{equation*}
\end{corollary}
From \eqref{18}, we note that
\begin{equation}
T_{n+2}\big(x\big|a,b\big)+axT_{n+1}\big(x\big|a,b\big)+bT_{n}\big(x\big|a,b\big)=0.\label{19}
\end{equation}
Let
\begin{equation*}
T_{n}\big(x\big|a,b\big)=X.
\end{equation*}
Then, by \eqref{19}, we get
\begin{equation*}
X^{n+2}+axX^{n+1}+bX^{n}=0\quad \textrm{if and only if}
\end{equation*}
\begin{equation}
X^{2}+axX+b=0.\label{20}
\end{equation}
Thus, from \eqref{20} we have
\begin{equation}
X=\frac{-ax\pm\sqrt{(ax)^{2}-4b}}{2}.\label{21}
\end{equation}
Hence, by \eqref{21}, we get
\begin{equation}
T_{n}\big(x\big|a,b\big)=C_{1}\bigg(\frac{-ax+\sqrt{(ax)^{2}-4b}}{2}\bigg)^{n+1}+C_{2}\bigg(\frac{-ax-\sqrt{(ax)^{2}-4b}}{2}\bigg)^{n+1}.\label{22}
\end{equation}
From \eqref{18} and \eqref{22}, we note that
\begin{equation}
C_{1}=\frac{1}{\sqrt{(ax)^{2}-4b}}\quad\mathrm{and}\quad C_{2}=-\frac{1}{\sqrt{(ax)^{2}-4b}}.\label{23}
\end{equation}
Therefore, by \eqref{22} and \eqref{23}, we obtain the following theorem.
\begin{theorem}
For $n\ge 0$, we have
\begin{equation*}
T_{n}\big(x|a,b\big)=\frac{1}{\sqrt{(ax)^{2}-4b}}\bigg(\bigg(\frac{-ax+\sqrt{(ax)^{2}-4b}}{2}\bigg)^{n+1}-\bigg(\frac{-ax-\sqrt{(ax)^{2}-4b}}{2}\bigg)^{n+1}\bigg).
\end{equation*}
\end{theorem}
From \eqref{13}, we note that
\begin{align}
\sum_{n=0}^{\infty}K_{n}\big(x\big|a,b,c\big)t^{n}
&=\big(1+cxt\big)\sum_{m=0}^{\infty}(-1)^{m}t^{m}\sum_{k=0}^{\infty}\binom{m}{k}(ax)^{m-k}b^{k}t^{k} \label{24} \\
&=\big(1+cxt\big)\sum_{n=0}^{\infty}\sum_{k=0}^{[\frac{n}{2}]}\binom{n-k}{k}(-1)^{n-k}(ax)^{n-2k}b^{k} t^{n} \nonumber\\
&=\sum_{n=0}^{\infty}\sum_{k=0}^{[\frac{n}{2}]}\binom{n-k}{k}(-1)^{n-k}(ax)^{n-2k}b^{k}t^{n}\nonumber\\
&\quad+\sum_{n=1}^{\infty}cx\sum_{k=0}^{[\frac{n-1}{2}]}\binom{n-k-1}{k}(-1)^{n-k-1}a^{n-2k-1}x^{n-2k-1}b^{k}t^{n}\nonumber\\
&=1+\sum_{n=1}^{\infty}\bigg(\sum_{k=0}^{[\frac{n}{2}]}\binom{n-k}{k}(-1)^{n-k}(ax)^{n-2k}b^{k} \nonumber\\
&\qquad +c\sum_{k=0}^{[\frac{n-1}{2}]}\binom{n-k-1}{k}(-1)^{n-k-1}a^{n-2k-1}x^{n-2k}b^{k}\bigg)t^{n}. \nonumber
\end{align}
Therefore, by \eqref{24}, we obtain the following theorem.
\begin{theorem}
For $n\ge 1$, we have
\begin{equation*}
\begin{aligned}
K_{n}\big(x\big|a,b,c)&= \sum_{k=0}^{[\frac{n}{2}]}\binom{n-k}{k}(-1)^{n-k}(ax)^{n-2k}b^{k}\\
&\quad +\frac{c}{a}\sum_{k=0}^{[\frac{n-1}{2}]}\binom{n-k-1}{k}(-1)^{n-k-1}(ax)^{n-2k}b^{k}.
\end{aligned}
\end{equation*}
\end{theorem}
Note that
\begin{equation}
T_{2n}\big(x\big|1,-1\big)=\sum_{k=0}^{n}\binom{2n-k}{k}x^{2(n-k)}=\sum_{r=0}^{n}\binom{n+r}{n-r}x^{2r}. \label{25}	
\end{equation}
The Morgan-Voyce polynomials are given by \cite{13, 14, 17, 18}
\begin{equation}
b_{n}(x)=xB_{n-1}(x)+b_{n-1}(x),\quad B_{n}(x)=(x+1)B_{n-1}(x)+b_{n-1}(x),\ (n\ge 1), \label{26}
\end{equation}
and
\begin{equation*}
b_{0}(x)=B_{0}(x)=1.
\end{equation*}
From \eqref{26}, we have
\begin{equation}
B_{n}(x)=\sum_{k=0}^{n}\binom{n+k+1}{n-k}x^{k},\quad b_{n}(x)=\sum_{k=0}^{n}\binom{n+k}{n-k}x^{k}.\label{27}
\end{equation}
Thus, by \eqref{25} and \eqref{27}, we get
\begin{equation}
T_{2n}\big(x\big|1,-1\big)=b_{n}(x^{2}),\ (n\ge 0).\label{28}
\end{equation}
By \eqref{26} and \eqref{28}, we get
\begin{equation}
x^{2}B_{n}(x^{2})=b_{n+1}(x^{2})-b_{n}(x^{2})=T_{2n+2}\big(x\big|1,-1\big)-T_{2n}\big(x\big|1,-1\big). \label{29}	
\end{equation}
Therefore, by \eqref{29}, we obtain the following theorem.
\begin{theorem}
For $n\ge 0$, we have
\begin{equation*}
x^{2}B_{n}(x^{2})=T_{2n+2}\big(x\big|1,-1\big)-T_{2n}\big(x\big|1,-1\big).
\end{equation*}
\end{theorem}
From \eqref{17}, we note that
\begin{align}
t\sum_{n=0}^{\infty}T_{n}\big(x\big|a,b\big)t^{n}&=\frac{t}{1+axt+bt^{2}}=t\sum_{m=0}^{\infty}(-1)^{m}t^{m}(ax+bt)^{m}\label{30} \\
&=t\sum_{m=0}^{\infty}(-1)^{m}t^{m}\sum_{k=0}^{\infty}\binom{m}{k}(ax)^{m-k}b^{k}t^{k} \nonumber\\
&=t\sum_{n=0}^{\infty}\sum_{k=0}^{[\frac{n}{2}]}\binom{n-k}{k}(ax)^{n-2k}(-1)^{n-k}b^{k}t^{n}\nonumber\\
&=\sum_{n=1}^{\infty}\sum_{k=0}^{[\frac{n-1}{2}]}\binom{n-k-1}{k}(ax)^{n-2k-1}(-1)^{n-k-1}b^{k}t^{n}.\nonumber
\end{align}
Thus, by \eqref{30}, we get
\begin{align}
&\sum_{n=1}^{\infty}\sum_{k=0}^{[\frac{n-1}{2}]}\binom{n-k-1}{k}(ax)^{n-2k-1}(-1)^{n-k-1}b^{k}t^{n}&=t\sum_{n=0}^{\infty}T_{n}\big(x\big|a,b\big)t^{n}\label{31}\\
&=\sum_{n=1}^{\infty}T_{n-1}\big(x\big|a,b\big)t^{n}. \nonumber
\end{align}
Therefore, by \eqref{31}, we obtain the following theorem.
\begin{theorem}
For $n\ge 1$, we have
\begin{equation}
T_{n-1}\big(x\big|a,b\big)= \sum_{k=0}^{[\frac{n-1}{2}]}\binom{n-k-1}{k}(-1)^{n-k-1}(ax)^{n-2k-1}b^{k}.\label{32}
\end{equation}
\end{theorem}
From \eqref{4} and \eqref{32}, we have the following corollary.
\begin{corollary}
For $n\in\mathbb{N}$, we have
\begin{equation*}
T_{n-1}\big(x\big|-1,-1\big)=f_{n}(x).
\end{equation*}
In particular, if $x=1$, we get
\begin{equation*}
T_{n-1}\big(1\big|-1,-1\big)=f_{n}.
\end{equation*}
\end{corollary}
From \eqref{6}, \eqref{7} and \eqref{13}, we note that
\begin{align}
\sum_{m=0}^{\infty}K_{m}\big(x\big|-1,0,0\big)(e^{t}-1)^{m}&=\frac{1}{1-x(e^{t}-1)}=\sum_{k=0}^{\infty}x^{k}\big(e^{t}-1\big)^{k}\label{33}\\
&=\sum_{k=0}^{\infty}x^{k}k!\frac{1}{k!}\big(e^{t}-1\big)^{k}=\sum_{k=0}^{\infty}x^{k}k!\sum_{n=k}^{\infty}{n \brace k}\frac{t^{n}}{n!} \nonumber\\
&=\sum_{n=0}^{\infty}\sum_{k=0}^{n}x^{k}k!{n \brace k}\frac{t^{n}}{n!}=\sum_{n=0}^{\infty}F_{n}(x)\frac{t^{n}}{n!}. \nonumber
\end{align}
On the other hand, by the generating function of the Stirling numbers of the second kind, we get
\begin{align}
\sum_{m=0}^{\infty}K_{m}\big(x\big|-1,0,0\big)m!\frac{1}{m!}\big(e^{t}-1\big)^{m}&=\sum_{m=0}^{\infty}K_{m}\big(x\big|-1,0,0\big)m!\sum_{n=m}^{\infty}{n \brace m}\frac{t^{n}}{n!}	\label{34}\\
&=\sum_{n=0}^{\infty}\sum_{m=0}^{n}K_{m}\big(x\big|-1,0,0\big)m!{n \brace m}\frac{t^{n}}{n!}. \nonumber
\end{align}
Therefore, by \eqref{33} and \eqref{34}, we obtain the following theorem.
\begin{theorem}
For $n\ge 0$, we have
\begin{equation*}
F_{n}(x)=\sum_{m=0}^{n}K_{m}\big(x\big|-1,0,0\big){n \brace m}m!.
\end{equation*}
\end{theorem}

\section{The Euler-Seidel matrix}
The Euler-Seidel matrix $(a_{n,k})_{n,k\ge 0}$ is defined by \cite{5, 10}
\begin{equation}
a_{0,n}(x)=a_{n}(x),\ (n\ge 0),\quad a_{k,n}(x)=a_{k-1,n}(x)+a_{k-1,n+1}(x),\ (n\ge 0,\ k\ge 1).\label{35}
\end{equation}
From \eqref{35}, we have binomial identities given by
\begin{equation}
\begin{aligned}
&a_{n,0}(x)=\sum_{k=0}^{n}\binom{n}{k}a_{0,k}(x),\\
&a_{0,n}(x)=\sum_{k=0}^{n}\binom{n}{k}(-1)^{n-k}a_{k,0}(x).
\end{aligned}\label{36}
\end{equation}
We let
\begin{equation*}
A(t)=\sum_{n=0}^{\infty}a_{0,n}(x)t^{n},\quad \overline{A}(t)=\sum_{n=0}^{\infty}a_{n,0}(x)t^{n}.
\end{equation*}
Then \eqref{36} is equivalent to Euler's formula
\begin{equation*}
\overline{A}(t)=\frac{1}{1-t}A\bigg(\frac{t}{1-t}\bigg). 	
\end{equation*}
From \eqref{15} and \eqref{16}, we see that the Euler-Seidel matrix $\big(a_{ij}(x)\big)_{i,j\ge 0}$ associated with $\big(a_{0,n}(x)\big)_{n\ge 0}=\big(K_{n}(x|a,b,c)\big)_{n\ge 0}$ is given by
\begin{align*}
&\left(\begin{matrix}
1 & (c-a)x & (a^{2}-ac)x^{2}-b & \cdots \\
(c-a)x+1 & (a^{2}-ac)x^{2}+(c-a)x-b  & \vdots & \cdots \\
\big(a^{2}-ac\big)x^{2}+2(c-a)x-b+1 & \vdots & \vdots & \cdots \\
\vdots & \vdots & \vdots & \cdots
\end{matrix}\right).
\end{align*}
\emph{Remark.} The continued fractions are given by
\begin{equation*}
\cfrac{K_{n}\big(x\big|a,b,c\big)}{K_{n-1}\big(x\big|a,b,c\big)}=-ax+\cfrac{b}{ax+\cfrac{b}{ax+\cfrac{b}{ax+\cfrac{b}{ax+\cdots+\cfrac{b}{ax+\cfrac{b}{ax-cx} }   }  } }  }
\end{equation*}
and
\begin{equation*}
\cfrac{T_{n}\big(x\big|a,b\big)}{T_{n-1}\big(x\big|a,b\big)}=-ax+\cfrac{b}{ax+\cfrac{b}{ax+\cfrac{b}{ax+\cdots+\cfrac{b}{ax+\cfrac{b}{ax} } } }  }
\end{equation*}
\section{Conclusion}
In this work, we successfully introduced and characterized the generalized Chebyshev polynomials $K_{n}(x|a,b,c)$, establishing a single framework that unifies all four classical kinds of Chebyshev polynomials alongside related sequences. Through explicit closed-form expressions, tridiagonal determinant representations, and recurrence relations, we illuminated the underlying algebraic structures of these families. Furthermore, we demonstrated the broad scope of this generalization by linking it to Morgan-Voyce, Fibonacci, and Fubini polynomials, while analyzing its Euler-Seidel matrix properties and continued fraction representations. These findings offer a consolidated foundation for exploring further extensions and applications across combinatorics, special functions, and linear algebra.

\section{ CONFLICT OF INTERES}
On behalf of all authors, the corresponding author states that there is no conflict
of interest

\end{document}